\documentclass[12pt,reqno]{amsart}
\usepackage[final]{epsfig}
\newtheorem{theorem}{Theorem}[section]
\newtheorem{lemma}[theorem]{Lemma}

\newtheorem{corollary}[theorem]{Corollary}

\renewenvironment{proof}{\noindent {\bf Proof.}}{ \hfill\qed\\ }

\begin{document}

\def\e{\varepsilon}
\def\t{\theta}
\def\Q{{\bf Q}}
\def\Tt{T_{\theta}}
\def\rt{r_{\theta}}
\def\ft{f_{\theta}}
\def\St{S_{\theta}}
\def\Xt{X_{\theta}}
\def\wt{w_{\theta}}
\def\ct{c_{\theta}}
\def\Ct{C_{\theta}}
\def\x{\vec{x}}
\def\R{\mathbb{R}}
\def\N{\mathbb{N}}
\def\L{\mathcal{L}}
\def\A{\mathcal{A}}
\def\BL{\mathcal{BL}}
\def\P{\mathcal{P}}
\def\R{\mathbb{R}}
\def\em{\varepsilon_m}
\title{Complexity and growth for polygonal billiards}
\author{J.~Cassaigne, P.~Hubert and S.~Troubetzkoy}

\address{Institut de math\'ematiques de
Luminy, CNRS Luminy, Case 907,
F-13288 Marseille Cedex 9, France}
\email{cassaign@iml.univ-mrs.fr}
\urladdr{http://iml.univ-mrs.fr/{\lower.7ex\hbox{\~{}}}cassaign/}
\address{Institut de math\'ematiques de
Luminy, CNRS Luminy, Case 907,
F-13288 Marseille Cedex 9, France}
\email{hubert@iml.univ-mrs.fr}

\address{Centre de physique th\'eorique et Institut de math\'ematiques de
Luminy, CNRS Luminy, Case 907,
F-13288 Marseille Cedex 9, France}
\email{serge@cpt.univ-mrs.fr}
\email{troubetz@iml.univ-mrs.fr}
\urladdr{http://iml.univ-mrs.fr/{\lower.7ex\hbox{\~{}}}troubetz/}

\begin{abstract}
We establish a relationship between the word complexity and the number of 
generalized diagonals for a polygonal billiard.  We conclude that in
the rational case the complexity function has cubic upper and lower bounds.
In the tiling case the complexity has cubic asymptotic growth.
\end{abstract}
\maketitle
\section{Introduction}

A billiard ball, i.e.~a point mass,
moves inside a polygon $Q \subset \R^2$ with unit speed 
along a straight line until it reaches the boundary $\partial Q$, then 
instantaneously changes direction according to the mirror law: ``the angle of 
incidence is equal to the angle of reflection,'' and continues along the 
new line.

How complex is the game of billiards in a polygon? The first results
in this direction, proven independently by Sinai \cite{S} and Boldrighini,
Keane and Marchetti \cite{BKM} is that the metric entropy with respect
to the invariant phase volume is zero.  Sinai's proof in fact shows more,
the ``metric complexity'' grows at most polynomially.  Furthermore, it is known
that the topological entropy (in various senses) is zero \cite{K,GKT,GH}.

To prove finer results there are two natural quantities one can count,
one is the number of generalized diagonals, that is (oriented) 
orbit segments which begin and end in a vertex of the polygon and contain
no vertex of the polygon in their interior.  The number of links 
of a generalized diagonal is called its combinatorial length
while its geometric length is simply the
sum of the lengths of the segments.  
Let $N_g(t)$ (resp.~$N_c(n)$) be the number of generalized diagonals 
of geometric  (resp.~combinatorial) length at most $t$ (resp.~$n$).
Katok has shown that $N_g(t)$ grows slower than any
exponential \cite{K}.  Masur has shown that for rational polygons, that
is for polygons all of whose inner angles are commensurable with $\pi$,
$N_g(t)$ has quadratic upper and lower bounds \cite{M1,M2}.
By elementary reasoning there is a constant $B>1$ such that
$B^{-1} \le N_c(t)/N_g(t) \le B$, thus
all of these results easily extend to the quantity $N_c(n)$.  Furthermore,
Veech has shown that there is a special class of polygons now commonly
referred to as Veech polygons, for
example regular polygons, such that
the quantity $N_g(t)/t^2$ admits a limit as $t$ tends
to infinity \cite{V,V1}.  

To introduce the second natural quantity which can be counted, 
label the sides of $Q$ by symbols from a finite alphabet $\A$ whose 
cardinality is equal to the number of sides of $Q$.
We code the orbit by the sequence of sides it hits.  Consider
the set $\L(n)$ of all words of length $n$ which arise
via this coding.  Let $p(n) = \# \L(n)$, this is called the complexity
function of the language $\L(\cdot).$
The only general results known about the complexity function is that
it grows slower than any exponential \cite{K} and at least quadratically
\cite{Tr}.  For billiards in a square the
complexity function has been explicitly calculated, albeit for a slightly
different coding (the alphabet consists of two symbols, one for vertical
sides one for horizontal sides) \cite{Mi,BP}. For this coding of the
square the collection of codes which appear are known as the
Sturmian sequences.  In fact it is not hard
to relate the complexity functions for the two different codings, the
relationship is $p_4(n) = 4p_2(n) -4$.  There are some related
results on the complexity when one restricts to certain initial conditions:
for rational polygons the ``directional complexity'' in each direction
is known explicitly \cite{H1}, while for general polygons there are polynomial
upper bounds for the directional complexity \cite{GT}.

There are several good surveys of billiards in polygons, in these surveys
one can find  more details about the definitions and
more precise statements of the above mentioned results. We refer
the reader to \cite{Gu1,Gu2,MT,T}.

Our main theorem shows that $p(n)$ and $N_c(n)$ are related.

\begin{theorem}
For any convex 
polygon 
$$p(n) = \sum_{j=0}^{n-1} N_c(j).$$
\label{theorem1}
\end{theorem}

Here we remark that $N_c(0)$ is the number of vertices of the polygon
while the sides of $Q$ are not counted as generalized diagonals.
Applying the above mentioned results of Masur's \cite{M1,M2} we 
immediately conclude

\begin{corollary}
If $Q$ is a rational convex 
polygon then there are positive constants $D_1,D_2$
such that
$$D_1 < p(n)/n^3 < D_2$$
for all $n \in \N \backslash \{0\}.$
\end{corollary}

Next we exhibit several examples where there are exact asymptotics. We
show

\begin{theorem}\label{thm::2}
If $Q$ is the square, the isosceles right triangle or the equilateral triangle
then
\begin{equation}\label{e1}
\lim_{n \to \infty} \frac{p(n)}{n^3}
\end{equation}
exists. The following table 
expresses the limit.

\begin{center}
\begin{tabular}{|l|c||l|c||l|c|}\hline
 & & & & &\\
Square  &   $\displaystyle\frac{4}{\pi^2}$ & 
$\displaystyle\left (\frac{\pi}2,\frac{\pi}4,\frac{\pi}4\right )$-triangle &  
$\displaystyle \frac{2}{3\pi^2}$ &
Equilateral triangle & $\displaystyle\frac{3}{4\pi^2}$ \\ 
 & & & & & \\ \hline
\end{tabular}
\end{center}
\end{theorem}

The proof of Theorem \ref{theorem1} is 
split into two parts.  The first part is combinatorial.
It uses the notion of bispecial words which was developed by 
Cassaigne \cite{C}.  The second part is geometric and uses a counting
argument based on Euler's formula.

{\bf Remark:}
it is known that for $n$ sufficiently large
the complexity of each aperiodic individual word is
$4(n+1)$ for the square, $3(n+2)$ for the equilateral triangle,
$4(n+2) $ for the isosceles right triangle and $6(n+2) $ for the half 
equilateral triangle \cite{H,H1,H2}. For the square the complexity
is four times larger than that of Sturmian sequences, thus the
fact that $p(n)$ is asymptotically four times the number of Sturmian words
of length $n$ is not surprising \cite{Mi,BP}.

Any infinite word of eventual complexity $3(n+2)$ 
whose language (i.e.~the collection of finite factors) 
is invariant under cyclic permutations of the letters arises as
the coding of a billiard trajectory in the equilateral triangle \cite{H}.
The third entry of the table in 
Theorem \ref{thm::2} gives the asymptotic growth rates of
the number of all such words.

Two interesting tiling cases remaining to evaluate the limit (\ref{e1}) are the
$\displaystyle \left (\frac{\pi}2,\frac{\pi}3,\frac{\pi}6\right )$-triangle
and the hexagon. In this triangular case
the methods developed for the other three cases allow us to conclude
that this limit exists and to calculate it explicitly.  Since
the combinatorics of this case is more complicated than the others,
we leave its explicit computation to the dedicated reader.

The hexagonal case remains open since the corresponding lattice
point counting problem seems not to have been investigated.

It would be interesting to know if the limit (\ref{e1}) exists in the
case of Veech polygons and also to exhibit cases when it does not exist.

\section{A combinatorial lemma}

Let $p(0) := 0$ and for any $n \ge 1$ let $s(n) := p(n+1) - p(n).$
For $u \in \L(n)$ let
\begin{eqnarray*}
m_l(u) & := & \#\{a \in \A: au \in \L(n+1)\}\\
m_r(u) & := & \#\{b \in \A: ub \in \L(n+1)\}\\
m_b(u) & := & \#\{(a,b) \in \A^2: aub \in \L(n+2)\}.
\end{eqnarray*}
We remark that all three of these quantities are larger than or equal
to one.
A word $u \in \L(n)$ is called left special if $m_l(u) > 1$, right special
if $m_r(u) > 1$ and bispecial if it is left and right special.
Let $\BL(n) := \{u \in \L(n): u \text{ is bispecial}\}.$
In this section we show that 

\begin{theorem}\label{thm::com}
For any polygon $Q$
$$s(n+1) - s(n) = 
\sum_{v \in \BL(n)} \Big (m_b(v) - m_l(v) - m_r(v) + 1 \Big ).$$
\label{combi}
\end{theorem}

{\bf Remark:} there is no assumption of convexity for this theorem, in fact
it is not necessary that the language arises from the coding of a 
polygonal billiard.

\begin{proof}
Since for every $u \in \L(n+1)$ there exists $b \in \A$ and
$v \in \L(n)$ such that $u = vb$ we have
$$s(n) = \sum_{u \in \L(n)} (m_r(u) - 1).$$
Thus
$$s(n + 1) - s(n)  = \sum_{v \in \L(n+1)} \Big (m_r(v) - 1 \Big ) - 
\sum_{u \in \L(n)} \Big (m_r(u) - 1 \Big ).$$
For $u \in \L(n+1)$ we can write $u = av$ where $a \in \A$ and 
$v \in \L(n)$, thus
$$s(n+1) - s(n) = \sum_{v \in \L(n)} \left [ \sum_{av \in \L(n+1)}
\Big (m_r(av) - 1 \Big ) - \Big (m_r(v) -1 \Big ) \right ].$$

For any word $v \in \L(n)$ and $av \in \L(n+1)$ any legal prolongation
to the right of $av$ is a legal prolongation to the right of $v$ as well
thus if $m_r(v) = 1$ then $m_r(av) = 1$.  Thus words with $m_r(v) = 1$
do not contribute to the above sum.
Thus $s(n+1) - s(n)$ is equal to the above sum
restricted to those $v$ such that $m_r(v) > 1.$ 
If furthermore $m_l(v) = 1$ then there is only a single $a$ such
that $av \in \L(n+1)$. For this $a$ we have $m_r(av) = m_r(v)$ thus
such words do not contribute to the sum either.  Thus we can restrict
the sum to bispecial words, yielding

$$s(n+1) - s(n) = \sum_{v \in \BL(n)} 
\left [ \sum_{av \in \L(n+1)} \Big (m_r(av) - 1 \Big ) - \Big (m_r(v) -1 \Big ) \right ].$$

The lemma follows since for any $v \in \BL(n)$ we have
$$ m_b(v) = \sum_{av \in \L(n+1)} m_r(av)$$
and
$$ m_l(v) = \sum_{av \in \L(n+1)} 1.$$

\end{proof}

\section{Proof of theorem \ref{theorem1}}

Let $X := \{ (s,v): s \in \partial Q \text{ and } v \text{ is an inner
pointing unit vector} \}$ and $P$ the ``partition'' of $X$ induced by the
sides of $Q$. The ambiguity of the definition of $P$ at the vertices plays
no role in our discussion.
Let $T: X \to X$ be the billiard ball map.  An element of the partition
$P \vee T^{-1}P \vee \cdots \vee T^{-n+1} P$ is called an $n$--cell.
The code of every point in an $n$--cell has the same prefix of length $n$, 
thus there is a bijection between the set of $n$--cells and the language
$\L(n)$.

If the footpoint of $T^nx$ is a vertex 
then we say that $x$ belongs to a discontinuity of order $n$.
A discontinuity (of any order) is locally a curve whose endpoints
lie on the boundary of $X$ or on a discontinuity of lower order.
We call each piece between such endpoints a smooth branch of the
discontinuity.

\begin{figure}
\centerline{\psfig{file=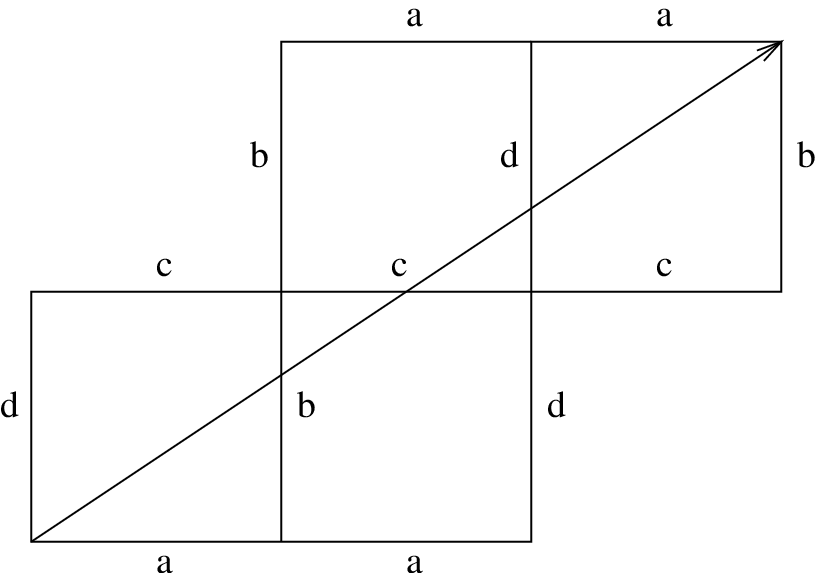,height=50mm}}
\caption{A generalized diagonal of combinatorial length 4 with code
$bcd$.}
\end{figure}    

For $v \in \BL(n)$ let $gd(v)$ be the number of generalized diagonals
of length $n+1$ such that the code of (the nonsingular part of)
the generalized diagonal is $v$
(see figure 1). Let $I_l(v) := m_l(v) - 1$ and $I_r(v) := m_r(v) - 1$.
For short we call  $(I_l(v),I_r(v),gd(v))$ the index of $v$.

\begin{lemma}\label{lem::geo} Suppose that $Q$ is a convex 
polygon.
For any $v \in \BL(n)$ 
$$m_b(v) = I_l(v) + I_r(v) + gd(v) +1$$
\end{lemma}

\begin{proof}
We consider the $n$--cell $C$ with bispecial code $v$. 
Note that an $n$--cell is a convex polygon \cite{K}, thus
geometrically the number
$m_b(v)$ corresponds to the number of pieces $C$ is cut into by the
discontinuities of order $-1$ and $n$.
  
Let $r$ be the number of sides of $Q$.
There are $I_l(v) \le r$ vertices of $Q$ which produce the splitting on 
the left,
they cut $C$ via singularities of $T^{-1}$.  Similarly there are
$I_r(v) \le r$ vertices which produce the splitting on the right, they
correspond to cutting $C$ via singularities of $T^{n}$.

Suppose the index of $v$ is $(i,j,k)$.  The cell $C$ is cut by $i+j$
singularities with $k$ intersections inside the interior of $C$.  
We claim that since $Q$ is convex, each of these $k$ intersections consists 
of  an intersection of exactly two smooth branches of the singularities.
Consider an intersection point $x$. Its forward orbit arrives at a vertex 
in say $m > 0$ steps and ends.  Thus $x$ belongs to the interior of a 
discontinuity of order $m$. The forward orbit hits no other vertex before 
time $m$, and by definition ends at time $m$, thus $x$ belongs to the 
interior of no other discontinuity of positive order.  There are two
possible continuations by continuity of the orbit of $x$.  If either
of these continuations is a generalized diagonal or tangent to a side
of $Q$ then $x$ is an end point of another singularity of positive order.
In the second case the order of this additional singularity is also $m$,
while in the first case it is strictly larger than $m$.
We note that the second possibility can only happen if $Q$ is not convex.
Similarly, considering the backwards orbit we see that $x$ belongs to 
the interior of a single discontinuity of negative order. If $Q$ is not 
convex then it is not the end point of any negative discontinuity of greater 
order. The claim is proven.

Next we will use Euler's formula to conclude our lemma.  
Let $F,E,V$ stand for the number of
faces, edges and vertices respectively of the partition of the interior
of $C$ by the discontinuities of order $-1$ and $n$.  
We have $E := i + j + 2k$ and $V:=k$.  
By Euler's formula we have
$V-E+F = 1$ thus $F = 1-V+E= 1+i+j+k$. As discussed above $m_b(v) = F$.
\end{proof}

\noindent
{\bf Proof of Theorem \ref{theorem1}.}
The theorem follows immediately from lemma \ref{lem::geo} and Theorem
\ref{thm::com}
since $N_c(n) = \sum_{j=0}^{n-1} \sum_{v \in \BL(j)} gd(v)$.
\hfill\qed

\section{Proof of Theorem \ref{thm::2}}

It is well known that if the images of $Q$ under the
action $A(Q)$ tile the plane, then $Q$ is the square, the equilateral triangle,
the right isosceles triangle or the half equilateral triangle (i.e.~the
triangle with angles $(\pi/2,\pi/3,\pi/6))$.  We will use this tiling 
to calculate $N_c(n)$.

\subsection{The square}
The tiling is the usual square grid.  
Fix a corner of the square and call it the origin of the grid. 
Consider all the generalized diagonals in the grid of combinatorial
length at most $n$ which
{\bf start} from this corner and are in the first quadrant. 
From figure 2 it is clear that 
the number $M_c(n)$ of such generalized diagonals is 
$$\#\left \{(i,j) \in \mathbb{N}^2: i+j \le n+1 \text{ and } 
\langle i,j \rangle = 1\right \}$$
where $\langle i,j \rangle$ is the gcd of $i$ and $j$. 
The condition $\langle i,j \rangle = 1$ arises since
generalized diagonals
stop as soon as they hit a vertex, thus if a line through the origin
hits several vertices (for example the line $y=x$),
it corresponds to only one generalized diagonal (starting at the origin).
Thus we only count it once.
Since there are four possible starting corners we have $N_c(n) = 4 M_c(n)$.

\begin{figure}
\centerline{\psfig{file=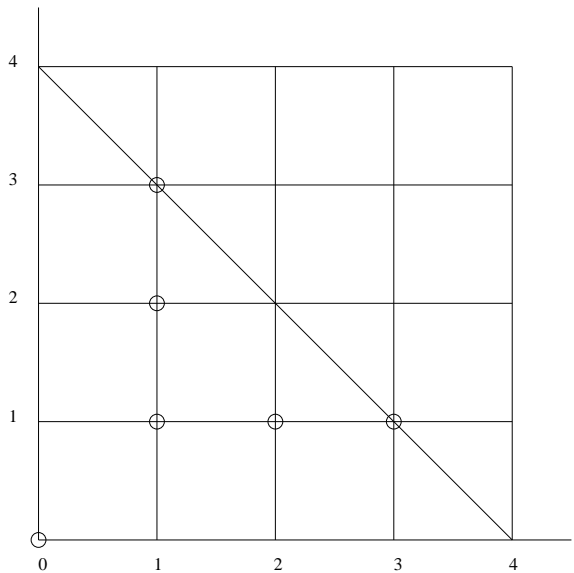,height=75mm}}
\caption{Counting the generalized diagonals starting at the origin
of combinatorial length at most 3 in the square.}
\end{figure}

The asymptotics of this quantity is well known by Mertens formula
and its generalizations \cite{HW,N}:

$$N_c(n) \sim \frac{12}{\pi^2}n^2.$$
Applying Theorem \ref{theorem1} we have

$$p(n) \sim \frac{12}{\pi^2} \sum_{k=1}^n k^2 \sim \frac{4}{\pi^2}n^3.$$

\subsection{The equilateral triangle}

We consider the images of $Q$ under the action of $A(G)$ specifying that
one of the vertices is at the origin and another at the point $(1,0)$.
We transform this grid to the grid in figure 3
via the affine mapping which fixes the vector 
{\Large \lower2pt\hbox{(}\hspace{-4pt}\tiny 
\begin{tabular}{c}
1 \\ 0
\end{tabular}\Large\hspace{-4pt}\lower2pt\hbox{)}}
and takes the vector
{\Large \lower2pt\hbox{(}\hspace{-4pt}\tiny 
\begin{tabular}{c}
$\cos(\pi/3)$ \\ $\sin(\pi/3)$
\end{tabular}\Large\hspace{-4pt}\lower2pt\hbox{)}}
to
{\Large \lower2pt\hbox{(}\hspace{-4pt}\tiny 
\begin{tabular}{c}
0 \\ 1
\end{tabular}\Large\hspace{-4pt}\lower2pt\hbox{)}}.

\begin{figure}
\centerline{\psfig{file=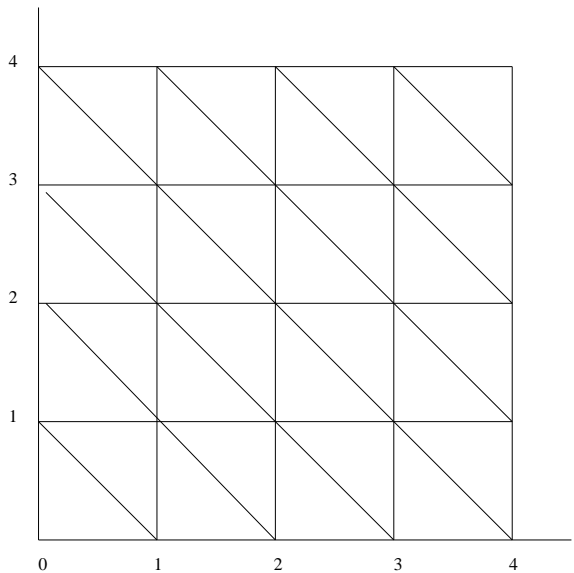,height=75mm}}
\caption{The affinely transformed grid of the equilateral triangle.}
\end{figure}

Consider all the generalized diagonals of combinatorial
length at most $n$ which
start from the origin and are in the first quadrant. 
Let  $M_c(n)$ be the cardinality of this set.  Since there
are 3 vertices we have
$N_c(n) = 3 M_c(n)$.  From figure 3 one sees that
$M_c(2n) = M_c(2n+1)$ since in traversing a square in the
tiling one alway crosses exactly two consecutive
copies of the fundamental triangle. From the figure it is also
clear that 
$M_c(2n) =\#\left \{(i,j) \in \mathbb{N}^2: i+j \le n+1 \text{ and } 
\langle i,j \rangle =1
\right\}.$ 

By Mertens formula \cite{HW,N}:

$$N_c(2n) = N_c(2n+1) \sim \frac{9}{\pi^2}n^2.$$
Applying Theorem \ref{theorem1} we have

$$p(n) \sim  \frac{3}{4\pi^2}n^3.$$

\subsection{The right isosceles triangle}
There are two different quantities which we must count.  First we
consider all the generalized diagonals of combinatorial
length at most $n$ which
start from the origin of the grid in figure 4a
and are in the first quadrant. 
Let  $M_1(n)$ be the cardinality of this set.  
We also
consider all the generalized diagonals of combinatorial
length at most $n$ which
start from the origin of the grid in figure 4b 
and are in the first octant. 
Let  $M_2(n)$ be the cardinality of this set.  There are two
vertices of our triangle with angle $\pi/4$ thus
$N_c(n) = M_1(n) + 2 M_2(n)$. 

\begin{figure}
\centerline{\psfig{file=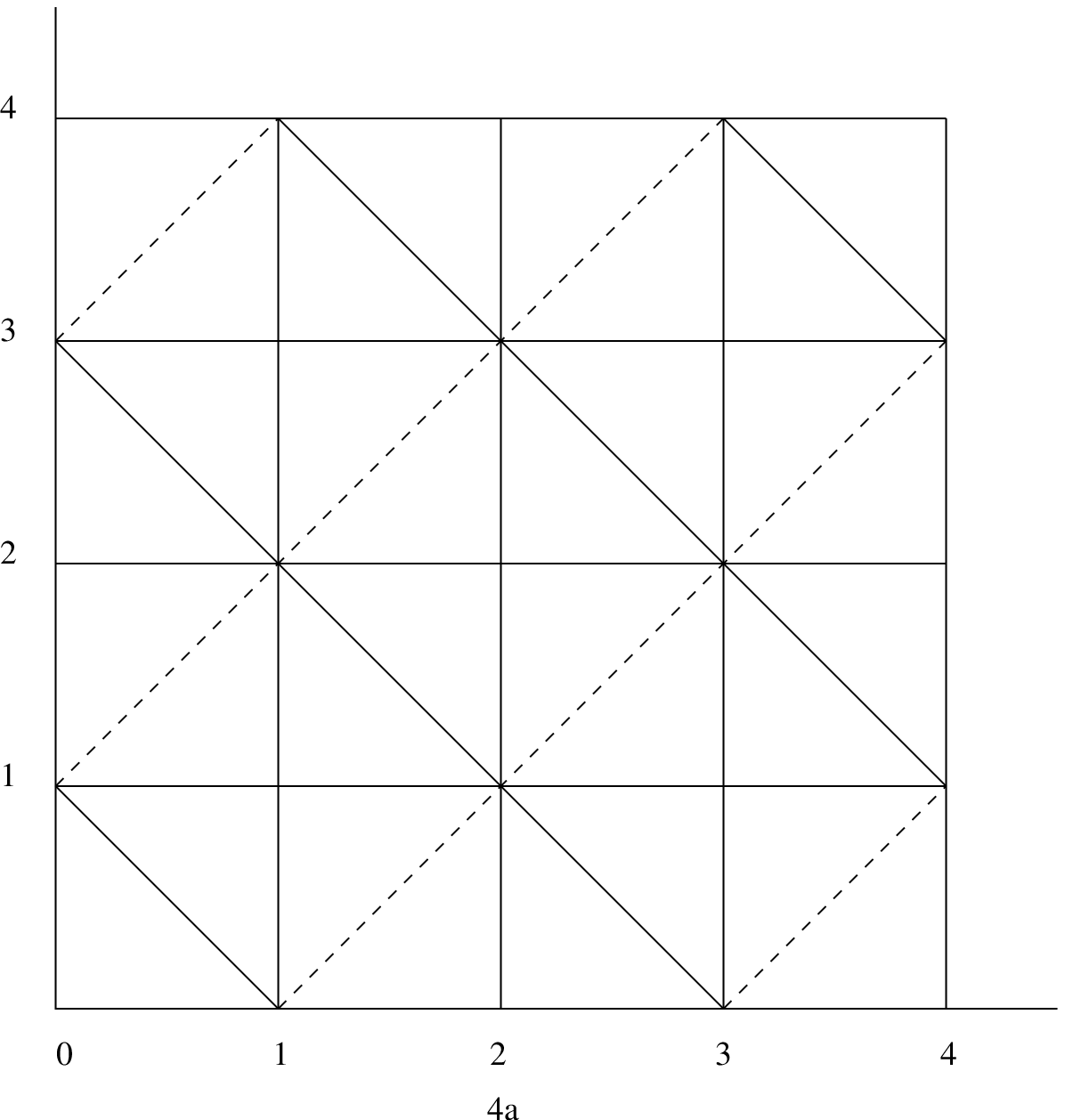,height=75mm} \quad \psfig{file=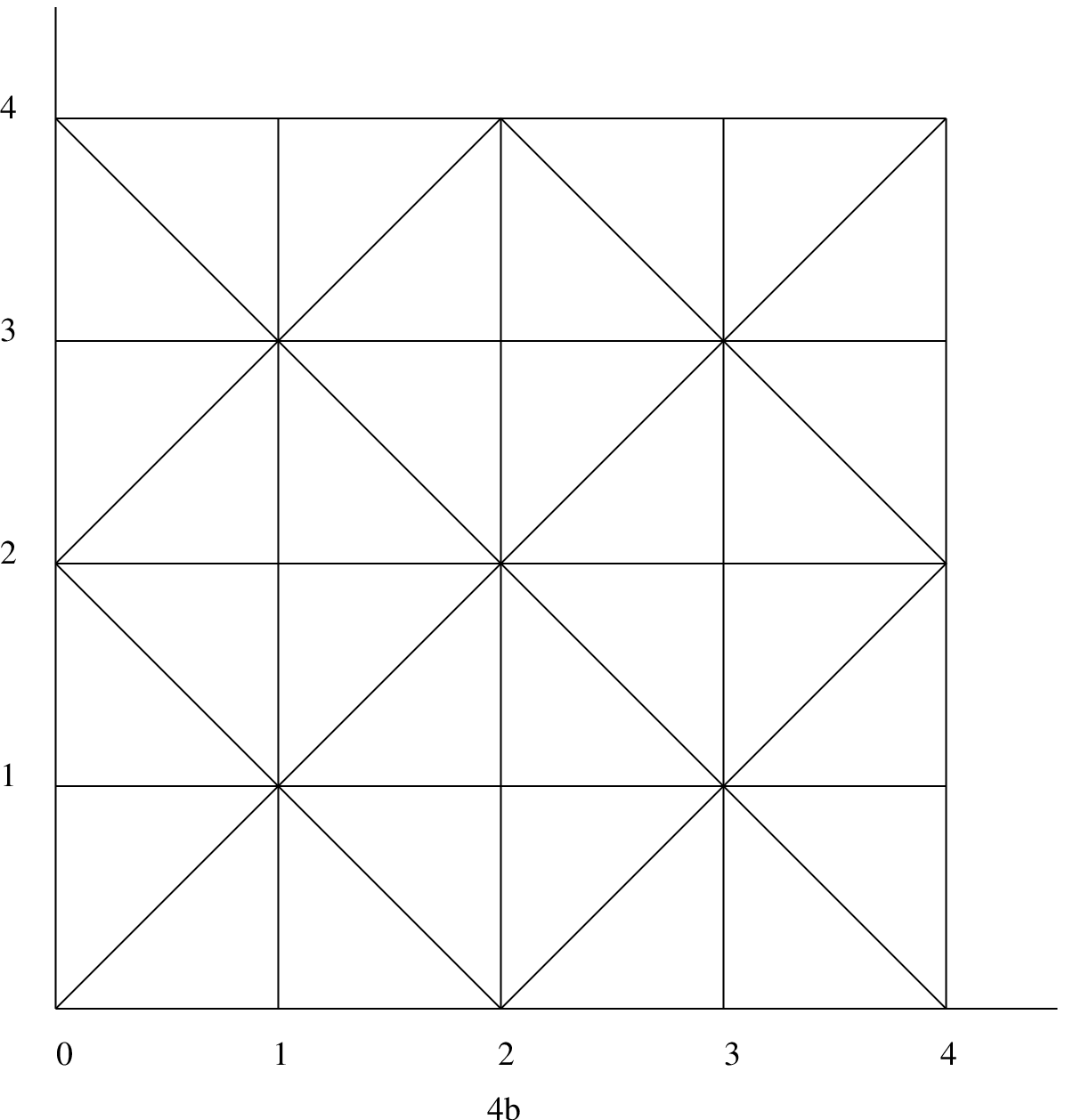,height=75mm}}
\caption{The grid of the right isosceles triangle.}
\end{figure}

The number $2 M_2(n)$ can also be interpreted as the cardinality of the set
of all the generalized diagonals of combinatorial
length at most $n$ which start from the origin of the grid in figure 4b 
and are in the first quadrant.
With this interpretation, if we overlay the grids from figures 4a and 4b 
we see that each generalized
diagonal counted in $M_1(n)$ is also a generalized diagonal counted in
$2 M_2(n+1)$ and each generalized diagonal counted in $2 M_2(n)$ is counted
in $M_1(n+1).$  Thus the asymptotics of $N_n(c)$ is the same as the
asymptotics of $2M_1(n)$.  

We want to count $M_1(n)$.  All generalized diagonals in the argument
below will start at the origin of the grid pictured in figure 4a and
all lengths will be combinatorial lengths.
We count first the solid lines which the generalized
diagonal crosses, we will deal later with the dashed lines. 
For most of the argument it will not matter whether the generalized diagonal
is simple or not (i.e.~contains no vertices in its interior), we
will restrict to the set of simple generalized diagonals only in the last
step of the proof.

Let $l(i,j)$ be the true combinatorial 
length of the generalized diagonal starting at the origin with end point
$(i,j)$ for any $(i,j) \in \mathbb{N}^2$. 
We view this length as the sum of the solid lines and the dashed lines
it crosses plus one.  If $i > j$ then the number of dashed lines it crosses is 
$\lfloor\frac{i-j}2\rfloor$. 

On the other hand the number of solid lines 
it crosses is characterized by the following statements.
Suppose that $n=3k$, then if it crosses $n-1$ solid lines then
$i + j = 2n/3 + 1 = 2k + 1$. Inversely, supposing that $i + j = 2k +1$,
then it crosses $n - 1 = 3k - 1$ solid lines.

Combining these two facts we have if
\begin{equation}\label{eq1}
(i,j) \in \mathbb{N}^2 \text{ and } i > j \text{ and } i + j = 2k+1
\end{equation}
then
$$l(i,j) = n + \left \lfloor \frac{i-j}2 \right \rfloor.$$

We need to calculate the region $\mathcal{R}(n)$ consisting
of all $(i,j) \in \mathbb{N}^2$ such that $i > j$
and $l(i,j) \le n.$
To do this fix $(i,j)$ as in (\ref{eq1}) and a natural number
$m \le i-1$.  We compare how many fewer dashed lines are crossed
by the generalized diagonal ending at $(i-m,j)$ than by
the generalized diagonal ending at $(i,j)$. This comparison yields
$l(i-m,j) = n +  \left \lfloor \frac{i-j}2 \right \rfloor - 2m
+ \em$ where $\em = m \mod 2$.
Thus $l(i-m,j) \le n$ if and only if $n +  \left \lfloor \frac{i-j}2 
\right \rfloor - 2m + \em \le n$.  
A simple computation yields the following two implications
\begin{eqnarray*}
m > \frac{i-j}4   \quad & \Longrightarrow & \quad l(i-m,j) \le n\\
m \le \frac{i-j}4 - \frac12 \quad & \Longrightarrow & \quad l(i-m,j) \ge n+1
\end{eqnarray*}

Let $m_0 := \min \{ m: \ l(i-m,j) \le n\}$.
From the above implications we have
$m_0 \le \frac{i-j}4 + 1$ and $m_0 \ge \frac{i-j}4 - \frac12$.
Let $\mathcal{D}(n)$ be the line $x = -y/2 + n/2$.  This line is the
``ideal boundary'' of the region $\mathcal{R}(n)$.  The
following computation shows that 
the distance of the true boundary from the ideal boundary is uniformly bounded:
$$d \Big ((i-m_0,j),\mathcal{D}(n) \Big ) \le d\left ((i-m_0,j),(-\frac j2 + \frac n2,j) \right) =
\left |i + \frac j2 - \frac n2 - m_0 \right | \le \frac54.$$

Let $\Delta^+(n)$ be the triangle whose boundaries are the $x$--axis,
the line $y=x$ and
the line $\mathcal{D}(n)$. By symmetry we also define a
region $\Delta^-(n)$ in the second octant (i.e.~we consider $i < j$).  
Let $\tilde{M}_1(n)$ be the number of simple generalized diagonals
starting at the origin whose other end is in the region $\Delta(n) := 
\Delta^+(n)
\cup \Delta^-(n)$.
Since the distance of the set $\Delta(n)$ from
the set $\mathcal{R}$ is uniformly bounded (in $n$) the asymptotics of
$\tilde{M}_1(n)$ and $M_1(n)$ are the same.
By symmetry $\tilde{M}_1(n)$ is twice the number of relatively prime
lattice points in the region $\Delta^+(n)$.
 
The area of $\Delta^+(n)$ is $n^2/12$. Thus
applying Mertens \cite{HW,N} yields
$$M_1(n) \sim \tilde{M}_1(n) \sim 2 \times \frac{n^2}{2\pi^2}.$$
Thus 
$$N_c(n) = 2 M_1(n) \sim \frac{2n^2}{\pi^2}$$ 
and applying Theorem \ref{theorem1} gives
$$p(n) \sim \frac{2}{3\pi^2}n^3.$$

\subsection{The half equilateral triangle}

The procedure is along the same lines as the previous examples.
We consider the affinely transformed grid 
similarly to the case of the equilateral triangle.
Counting the generalized diagonals which start at the origin 
reduces to an application of Merten's formula.  The explicit 
description of the region to
which Merten's formula must be applied is more complicated than
in the previous examples, thus we do not carry it out.

{\bf Acknowledgements:} We would like to thank Samuel Leli\`evre
a critical reading of an earlier version of this article.

\end{document}